\documentclass{elsart}

\usepackage{amssymb}
\usepackage{epsfig}
\usepackage[english]{babel}
\usepackage{bm,amsmath}
\usepackage{graphicx}

\begin{document}

\begin{frontmatter}



\title{Reversals of chance in paradoxical games}

\author[label1]{P. Amengual}
\ead{pau@imedea.uib.es}
\ead[url]{http://www.imedea.uib.es/$\sim$pau}
\author[label2]{P. Meurs}
\ead{pascal.meurs@uhasselt.be}
\author[label2]{B. Cleuren}
\author[label1]{R. Toral}
\address[label1]{Departament de F\'{\i}sica and Institut Mediterrani d'Estudis Avan\c{c}ats, Universitat de
 les Illes Balears, E07122 Palma de Mallorca, Spain.}
\address[label2]{Hasselt University, B-3590 Diepenbeek, Belgium.}

\begin{abstract}
We present two collective games with new paradoxical features when
they are combined. Besides reproducing the so--called Parrondo
effect, where a winning game is obtained from the alternation of
two fair games, a new effect appears, i.e., there exists a current
inversion when varying the mixing probability between the games.
We present a detailed study by means of a discrete--time Markov
chain analysis, obtaining analytical expressions for the
stationary probabilities for a finite number of players. We also
provide some qualitatively insight into this new current inversion
effect.
\end{abstract}

\end{frontmatter}

\section{Introduction}
In the past few years there has been an increasing interest in
what is known in the literature as Parrondo's paradox
\cite{ha99.2,ha02.1}. This phenomenon shows that the alternation
of two fair (or even losing) games can result in a winning game.
These so-called Parrondo games were originally defined as follows:
Game A is a simple coin tossing game, where the player wins or
loses one coin with probabilities $p^A$ and $1-p^A$ respectively.
For game B the winning probability depends on the capital of the
player modulo three, governed by the set of probabilities
$\{p_B^1,p_B^2,p_B^3\}$. In the (stochastic) combination of these
games, either game $A$ or $B$ is played, with probabilities
$\gamma$ and $1-\gamma$ respectively. The games are said to be
fair/losing/winning when on average the player's capital
stabilizes/decreases/increases.

In the Parrondo games, to which we will further refer as the
original game $A$ and $B$, the following parameter values were
used: $p^A=\frac{1}{2}-\epsilon$, $p_B^1=\frac{1}{10}-\epsilon$,
$p_B^2=p_B^3=\frac{3}{4}-\epsilon$. When $\epsilon = 0$ both games
A and B played separately are fair, whereas if $\epsilon >0$ both
games turn out to be losing. The Parrondo effect appears when the
stochastic $(0<\gamma<1)$ or periodic combination of these
fair/losing games results in a winning game.

These games were first devised in 1996 by the Spanish physicist
Juan M. R. Parrondo, who presented them in unpublished form in
Torino, Italy~\cite{p96.1}. They served as a pedagogical
illustration of the flashing ratchet~\cite{r02.1}, where directed
motion is obtained from the random or periodic alternation of two
relaxation potentials acting on a Brownian particle, none of which
individually produce any net flux. Only recently a quantitative
relation has been established between the Brownian ratchet and
Parrondo's games \cite{tam03.1,tam03.2}.

Cooperative versions of the games, played by a set of $N$ players,
have also been studied. In \cite{t01.1}, a set of $N$ players are
arranged in a ring and each round a player is chosen randomly to
play either game A or B. The original game $A$ is combined with a new
game $B$, for which the winning probability depends on the state (winner/loser)
of the nearest neighbors of the selected player. A player is said
to be a winner (loser) if he has won (lost) his last game. In
\cite{t02.1} again a set of $N$ players is considered, but
for this case game A is replaced by a redistribution process where
a player is chosen randomly to give away one coin of his capital
to another player. When combining this new game with the original
game B, the paradox is reproduced.

In this work we present a new version of collective games, where
besides obtaining the desired result of a winning game out of two
fair games, a new feature appears. The games show under certain
circumstances a current inversion when varying $\gamma$, i.e. the
value of the mixing probability $\gamma$ determines whether you
end up with a winning or a losing game $A+B$. To our knowledge,
this effect is new in the literature on paradoxical games, and in
the related field of Brownian ratchets, as we will discuss in
Sec.\ref{parrondo_inversion}.

\section{The games}\label{collective_games}

The games  will be played by a set of $N$ players. In each round,
a player is selected randomly for playing. Then, with
probabilities $\gamma$ and $1-\gamma$ respectively game $A$ or $B$
is played. Game $A$ is the original game in which the selected
player wins or loses one coin with probability $p^A$ and $1-p^A$
respectively. The winning probabilities in game $B$ depend on the
collective state of all players. Again, as in \cite{t01.1}, a
player is said to be a winner or a loser when he has won or lost
respectively his last game. More precisely, the winning
probability can have three possible values, determined by the
actual number of winners $i$ within the total number of players
$N$, in the following way
\begin{equation}\label{probabilities_game_B}
p^B_i \equiv\textrm{probability to win in game $B$ }= \left\{
\begin{array}{cc}
p_B^1\hspace*{0.3truecm} \text{if} &i>\lceil\frac{2N}{3}\rceil,\vspace*{0.2truecm}\\
p_B^2\hspace*{0.3truecm} \text{if} &\lceil\frac{N}{3}\rceil\leq i\leq \lceil\frac{2N}{3}\rceil,\vspace*{0.2truecm}\\
p_B^3\hspace*{0.3truecm} \text{if} &i<\lceil\frac{N}{3}\rceil.
\end{array}\right.
\end{equation}
where the brackets $\lceil x \rceil$ denote the nearest integer to the number $x$.

\subsection{Analysis of the games}

\begin{figure}[!htb]
\centerline{\epsfig{figure=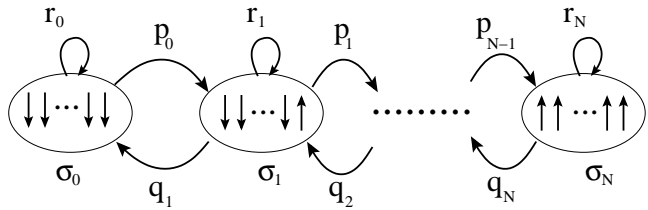,width=8.5cm}}
\caption{\label{states}Different states and allowed transitions
for $N$ players. The arrows indicate the state of each player
being a winner (\textit{arrow up}) or a loser (\textit{arrow
down}).}
\end{figure}

The main quantity of interest is the average gain of the
collection of $N$ players when playing the stochastic game $A+B$.
Since the winning probability of game $B$ only depend on the total
number of winners, it is sufficient to describe the games using a
set of $N+1$ different states
$\{\sigma_0,\sigma_1,\ldots,\sigma_N\}$. A state $\sigma_i$ is the
configuration where $i$ players are labeled as winner and $N-i$ as
loser. Transitions between the states will be determined by the
forward transition probability $p_i$, the backward transition
probability $q_i$, and the probability for remaining in the same
state $r_i$, see Fig.~\ref{states}.

Denoting as $P_i(t)$ the probability of finding the system in
state $\sigma_i$ at the t--{\textit{th}} round played, we can
write the equation governing its time evolution as
\begin{equation}\label{ME}
P_i(t+1)= p_{i-1} P_{i-1}(t) +r_i P_i(t) + q_{i+1} P_{i+1}(t),
\end{equation}
with $0\leq i\leq N$ and where the transition probabilities are
given by
\begin{eqnarray}
p_i & = &\frac{N-i}{N} \left[ \gamma~p^A+(1-\gamma)~p^B_i \right],\\
r_i & = &\frac{2i-N}{N} \left[ \gamma~p^A+(1-\gamma)~p^B_i \right]+\frac{N-i}{N},\\
q_i & = &\frac{i}{N} \left[
\gamma~(1-p^A)+(1-\gamma)~(1-p^B_i)\right].
\end{eqnarray}
The set of transition probabilities $(p_i,q_i,r_i)$ must satisfy
the normalization condition $p_i+r_i+q_i=1$, which implies for the
probabilities $P_i(t)$ that $\sum_{i=0}^{N}P_i(t)=1$, as long as
$\sum_{i=0}^{N}P_i(t=0)=1$.

This system of $N+1$ equations can be solved in the stationary
state, where the probabilities no longer depend on time
$P_i(t)=P^{\text{st}}_i$, and the general solution reads
\begin{eqnarray}\label{stat_prob}
P^{\text{st}}_i&=&\frac{1}{Z}p_0~p_1\cdots
p_{i-1}~q_{i+1}~q_{i+2}\cdots q_{N},
\end{eqnarray}
where $0\leq i\leq N$ and $Z$ is the normalization factor. Once
the stationary probabilities are calculated, we can obtain the
average winning probability over all states for the stochastic
combination $A+B$ (mixing probability $\gamma$) from
\begin{eqnarray}\label{pwin}
p_{\text{win}}^{A+B}=\sum_{i=0}^N \left[\gamma~
p^A+(1-\gamma)~p_i^B\right] P_i^{\text{st}}.
\end{eqnarray}
The average gain can then easily be evaluated through the
expression
$J^{A+B}=
2 p_{\text{win}}^{A+B}-1$.
The properties of the separate games A and B can be obtained by
replacing in the previous expressions $\gamma$ by $1$ or $0$ respectively.

\begin{table}
\centering
  \begin{tabular}{|c|c|}
  \hline
  \ \ \ $N$ \ \ \           &  $p_B^2$ \\   \hline\hline
   & \\
2 & $\frac{p_B^1-1}{p_B^1-p_B^3-1}$.        \\
   & \\
3 & $\frac{(p_B^1-1)(p_B^3+1)+\sqrt{(p_B^1-2)(p_B^1-1)p_B^3(p_B^3+1)}}{(p_B^1+p_B^3-1)}$\\
    & \\
4 &
$\frac{(p_B^1-1)^2(p_B^3+1)}{1+p_B^3+(p_B^1-2)(p_B^1+p_B^1p_B^3-(p_B^3)^2)}$\\
    & \\
5 & $\left[1-\frac{p_B^3}{p_B^1-1}\sqrt{\frac{5+2p_B^1(p_B^1-3)}{1+2p_B^3(1+p_B^3)}}\right]^{-1}$\\
    & \\
  \hline
  \end{tabular}
\caption{Condition on $p_B^2$ in order that game B is fair for
$N=2,\ldots,5$.}\label{tab:fair}
\end{table}

\subsubsection{The Parrondo effect}

\begin{figure}
\centerline{\epsfig{figure=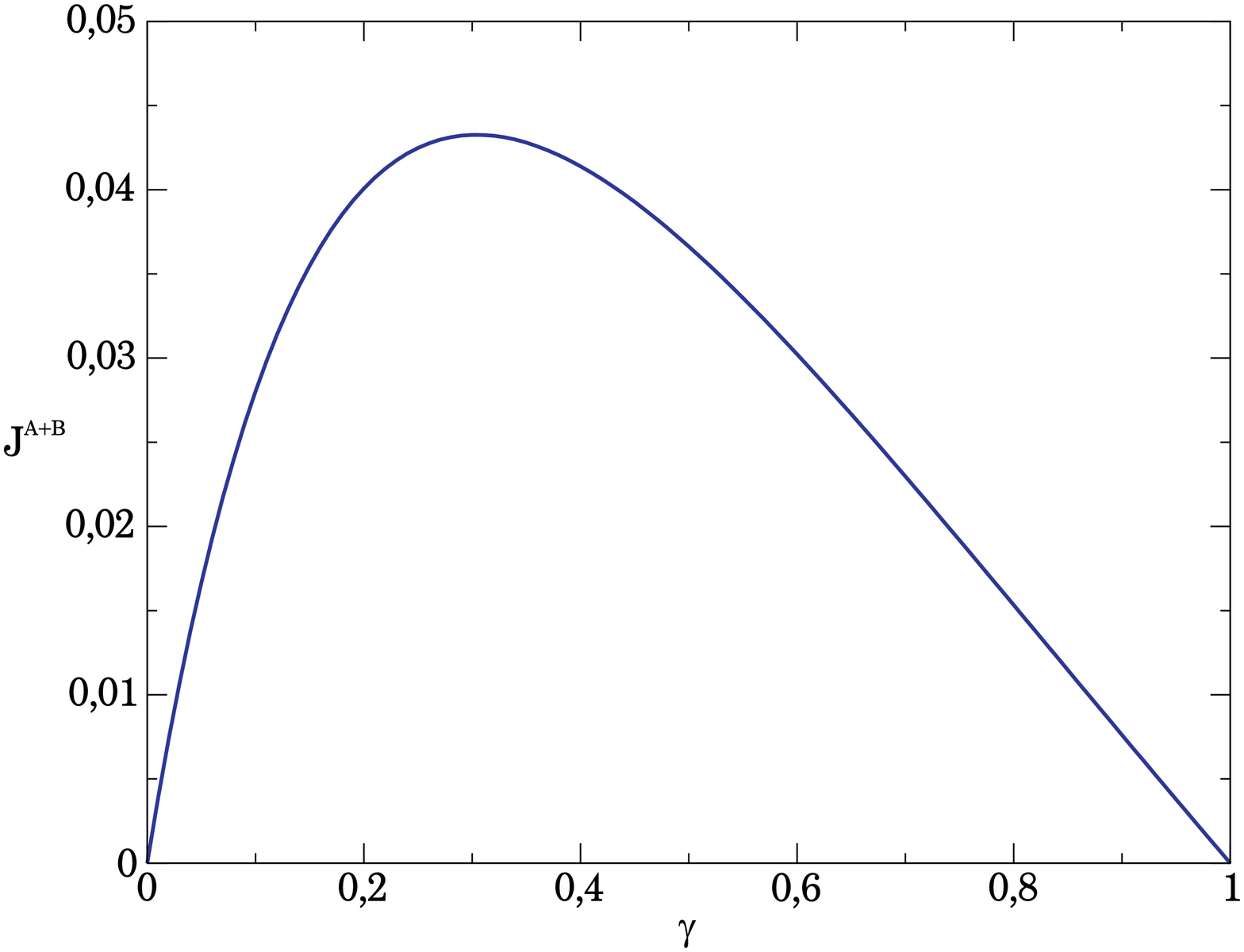,width=6.5cm}
\epsfig{figure=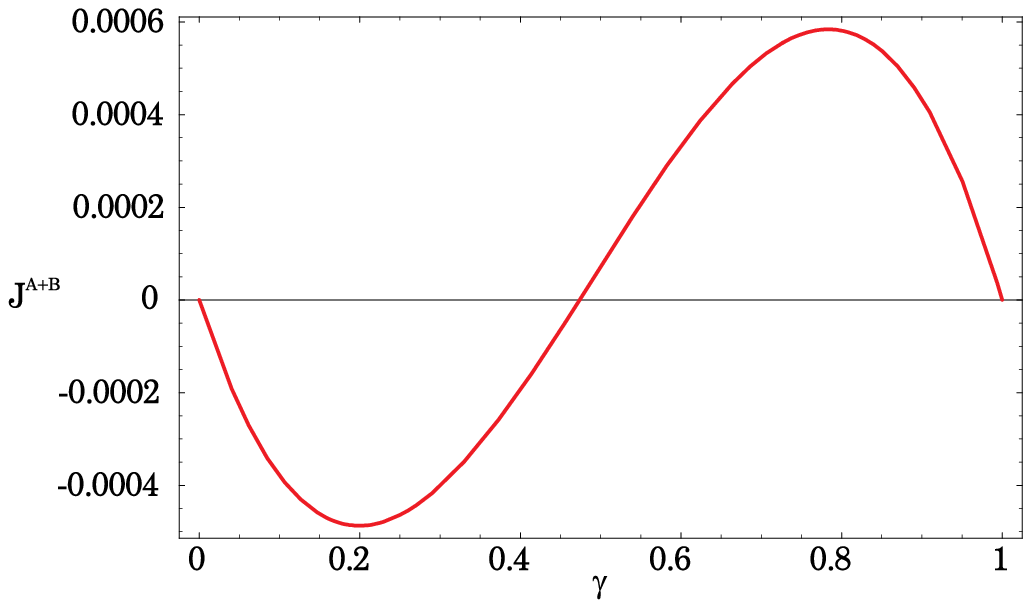,width=6.8cm,height =4.9cm }} \caption{\label{Parrondo_effect}\textit{\textbf{a)} }Plot
of the current versus the mixing probability $\gamma$ between
games $A$ and $B$ for $N=4$ with probabilities $p^A=\frac{1}{2}$,
$p_B^1=0.79$, $p_B^2=0.65$ and $p_B^3=0.15$. \textit{\textbf{b)}
}Plot of the current versus the mixing probability $\gamma$
between games A and B for $N=3$ with probabilities
$p^A=\frac{1}{2}$, $p_B^1=0.686$, $p_B^2=0.423$ and $p_B^3=0.8$. }
\end{figure}

The Parrondo effect appears when from the combination of two fair
games, we obtain a winning game. Clearly, game A is fair for
$p^A=1/2$. For game B the set of values $\{p_B^1,p_B^2,p_B^3\}$
giving a fair game is more difficult to determine because it
depends on the total number of players $N$. The conditions on
$p_B^2$ for a fair game B have been found analytically by a
symbolic manipulation program up to $N<13$. In Table
\ref{tab:fair} we find listed the conditions of fairness for
$p_B^2$ up to $N=5$.  When playing only game $B$ ($\gamma=0$), the
following symmetry in the stationary distribution can be deduced
from Eq.(\ref{stat_prob})

\begin{eqnarray}\label{symmetry}
P^{\text{st}, \{p_B^1,p_B^2,p_B^3\}}_i=P^{\text{st},
\{1-p_B^3,1-p_B^2,1-p_B^1\}}_{N-i}.
\end{eqnarray}

This property implies that $p_{\text{win}}^{A+B}$ is unaffected by
the parameter transformation:
$\{p_B^1,p_B^2,p_B^3\}\rightarrow\{1-p_B^3,1-p_B^2,1-p_B^1\}$. It
also means that for the parameter set
$\{p_B^1,p_B^2=1/2,1-p_B^1\}$, the stationary probability
distribution is symmetric over the states, i.e.
$P^{\text{st}}_i=P^{\text{st}}_{N-i}$. Therefore, when combining
this with game $A$, i.e., alternating two games with symmetric
probability distributions, always yields a fair game, independent
of the values of $\gamma$, $N$ and $p_B^1$. To see the Parrondo
effect, we need another, non-trivial, parameter set which yields a
fair game B.
For example, for $N=4$ we obtain a fair game B when
$p_B^1=0.79$, $p_B^2=0.65$ and $p_B^3 = 0.15$. The stochastic
combination with game $A$ reproduces the desired Parrondo effect,
see Fig.~\ref{Parrondo_effect}.a~.

\subsection{Results}

\subsubsection{Two players}

\begin{figure}[h]
\centerline{\includegraphics[width=6cm]{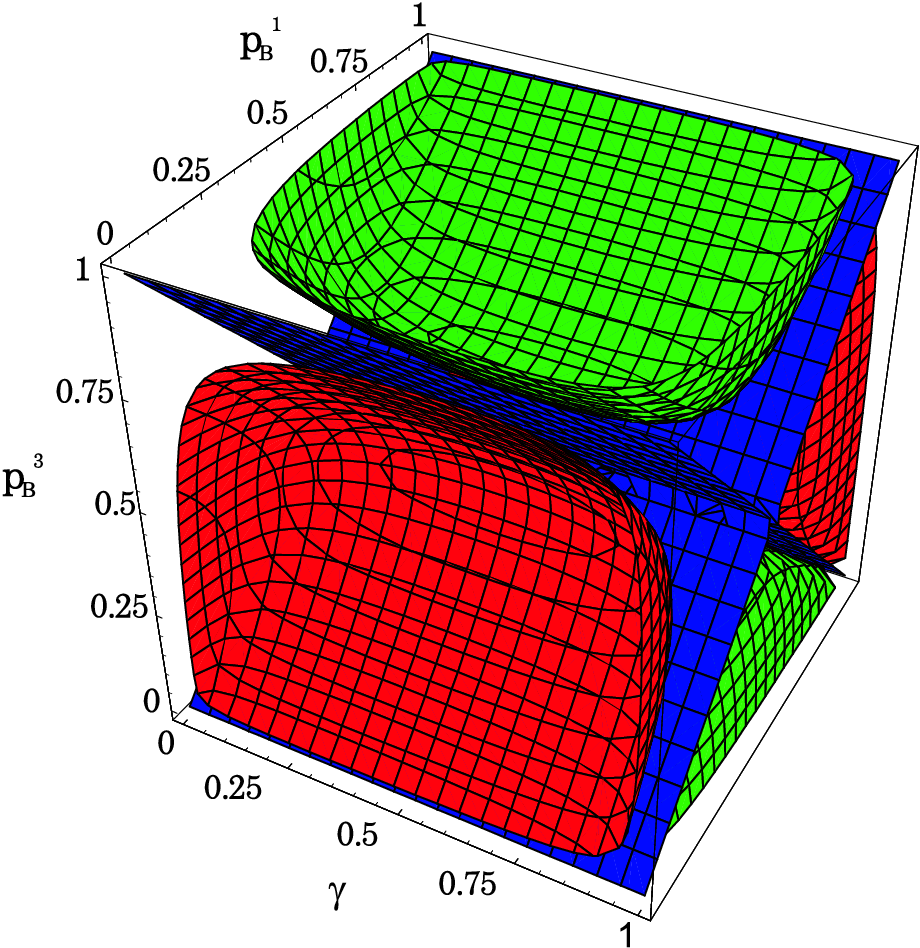}
\includegraphics[width=6cm]{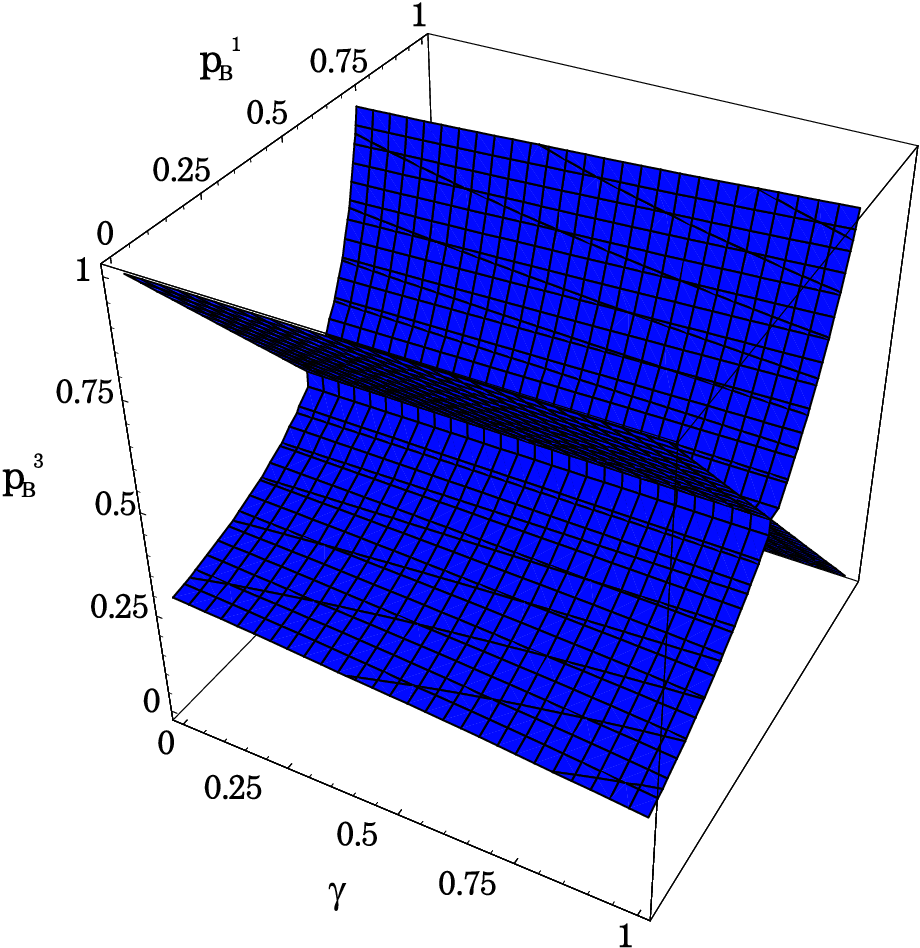}}
\caption{\label{N2_N3}\textit{\textbf{a)}} $N=2$. The regions in
parameter space for for which $p^{A+B}_{\textrm{win}}=0.5$,
$0.499$ and $0.501$, indicating the regions where $A+B$ is fair
(\textit{blue}), losing (\textit{red}) and winning
(\textit{green}) respectively. The blue diagonal planes show the
situations $p_B^1=1-p_B^3$ and $p_B^1=p_B^3$, for which $A+B$ is
fair, independent of $\gamma$. \textit{\textbf{b)}} $N=3$. The
regions in parameter space for which the mixing $(0<\gamma<1)$
between game $A$ and $B$ results in a fair game. Besides the
trivial diagonal plane, there is a curved plane -- not uniform in
$\gamma$ -- for which $J^{A+B}=0$.}
\end{figure}

For $N=2$ players, there are $3$ different states.
Fig.~\ref{N2_N3}.a shows the regions in parameter space
$\{\gamma,p_B^1,p_B^3\}$ where the mixing $(0<\gamma<1)$ between
games $A$ and $B$ results in a fair, winning or losing game. Note
that $p_B^2$ is fixed by the condition to have a fair game B, see
Table \ref{tab:fair}. Besides the case $p_B^1=1-p_B^3$, valid for
any number of players, also $p_B^1=p_B^3$ results in a fair game
for $N=2$, independent of the alternation probability $\gamma$.
From Eq.~(\ref{stat_prob}), one can deduce that $p_B^1=p_B^3$ and
$p_B^1=1-p_B^3$ imply a symmetric distribution $P^{\text{st}}_i$
over the states, i.e. $P^{\text{st}}_0=P^{\text{st}}_2$. As
mentioned before, this property prohibits any net current in the
system. For all other values of $p_B^1$ and $p_B^3$ the Parrondo
effect appears, that is, game $A+B$ is either a winning or a
losing game, cf. Fig.~\ref{N2_N3}.a.

\subsubsection{Three players}

Fig.~\ref{N2_N3}.b shows for $N=3$ the surfaces in parameter space
$\{\gamma,p_B^1,p_B^3\}$ where $A+B$ is a fair game. Besides the
plane $p_B^1=1-p_B^3$, there is a second, curved surface with
values of $\gamma$ different from $0$ and $1$ which results in
$J^{A+B}=0$. This curved surface is not uniform in $\gamma$ and is
therefore the collection of points of flux reversal between a
winning and losing game $A+B$. This implies that, depending on the
value of $\gamma$ we can either have a winning game or a losing
game by alternating between two fair games. For example, in
Fig.~\ref{Parrondo_effect}.b we have plotted the current $J^{A+B}$
vs. $\gamma$ for the set of probabilities $p^A=\frac{1}{2}$,
$p_B^1=0.686$, $p_B^2=0.423$ and $p_B^3=0.8$. For low values of
$\gamma$ the resulting game is a losing game, whereas for high
values of $\gamma$ the game turns to be a winning game, cf.
Fig~\ref{Parrondo_effect}.b. In both regions there exists an
optimal value for $\gamma$ giving a maximum current. We can
provide a qualitative picture that may help understanding the
mechanism by which the current inversion phenomenon takes place.

When playing exclusively game B ($\gamma=0$), the stationary
distribution $P^{\textrm{st}}_i$ is not homogeneous. This is
reflected by the fact that the central states
$\{\sigma_1,\sigma_2\}$ have a higher occupancy probability
($P^{st}_i$) than the boundary states $\{\sigma_0,\sigma_3\}$. On
the other hand, if we look to the winning probability, it is
higher in the latter set of states rather than in the former one
($p_B^1, p_B^3>p_B^2$).

Indeed, the central states can be labelled as \textit{losing}
states, as when combining game B with game A for any
$0\leq\gamma<1$, the average losing probability $p^l_i=\gamma
(1-p^A)+(1-\gamma)(1-p^B_i)<\frac{1}{2}$, i.e., it is more likely
on average for a player to lose money rather than to win when
being in one of these states. On the other hand, for the boundary
states the contrary is true: it is more likely to win money rather
than to lose for any $0\leq\gamma<1$, so we can refer to them as
\textit{winning} sites, i.e., $p^w_i=\gamma
p^A+(1-\gamma)p^B_i>\frac{1}{2}$.

When combining game B with A, the resulting game will be fair,
losing or winning depending on the net balance between the
occupancy probabilities and the average winning probability on
each set of central and boundary states. For low $\gamma$ values
(playing game $B$ more often), the high occupancy probability of
$\{\sigma_1,\sigma_2\}$ is the dominant part, and due to the low
winning probability on these sites the resulting game is a losing
game. On the contrary, for higher $\gamma$ values (playing game
$A$ more often), the winning probability on the boundary sites
$\{\sigma_0,\sigma_3\}$ is high enough to compensate their low
occupancy, resulting in a winning game.

\subsubsection{N players}

\begin{figure}
\centerline{\epsfig{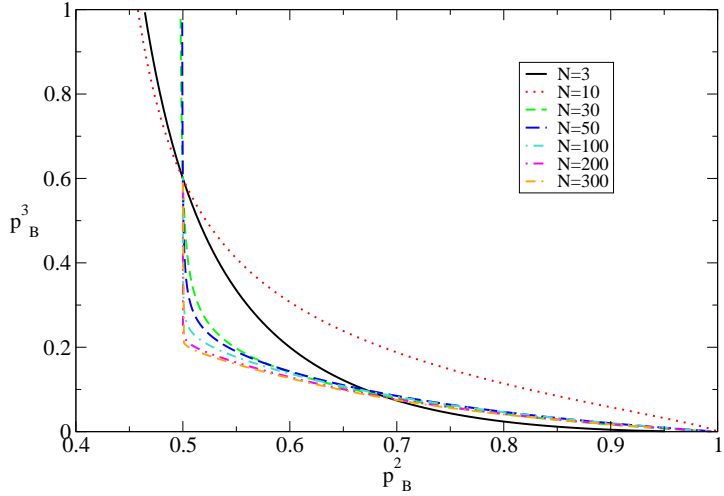}}
\caption{\label{parameter_spaces}Plot of the parameter space
$\{p_B^2,p_B^3\}$ for a fixed $p_B^1=0.4$ that gives a fair game B
for different values of $N =3$, $10$, $30$, $50$, $100$, $200$ and
$300$. As it can be seen, the curves seem to converge to a
limiting curve as $N$ increases.}
\end{figure}

For a general number of players, we have not been able to find the
analytical expressions for a fair game $B$. Nevertheless, we will
show numerically that the results for $N=3$ are representative for
any $N$. This is illustrated by Fig.~\ref{parameter_spaces}, where
the parameter space $\{p_B^2,p_B^3\}$ giving a fair game $B$ is
shown, corresponding to a fixed $p_B^1 = 0.4$ and different values
of $N$. As shown, the different curves seem to converge to a
limiting curve as $N$ increases. Note that all curves intersect at
the trivial point $\{p_B^1=0.4, p_B^2=0.5, p_B^3=0.6\}$.

\begin{figure}
\centerline{\epsfig{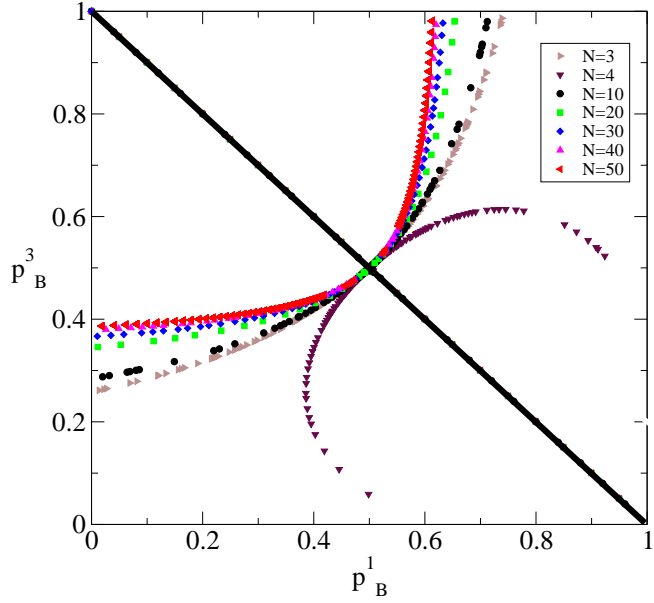}}
\caption{\label{fixed gamma}Plot of the points in parameter space
$\{p_B^1,p_B^3\}$ where (for $\gamma=0.4$ fixed) $A+B$ is a fair
game. Results for different values of the total number of players
$N = 3, 4, 10, 20, 30, 40$ and $50$ are shown. The diagonal line
shows the common plane $p_B^1=1-p_B^3$, that corresponds to a fair
game $B$ for any number of players $N$. }
\end{figure}

We can also obtain the parameter space where the current inversion
takes place, for different values of $N$. For clarity reasons we
show in Fig.~\ref{fixed gamma} only a vertical slice corresponding
to a fixed $\gamma =0.4$, and different values of $N$. Again, the
regions for which a flux inversion exists, doesn't seem to depend
much on $N$. The only exception is $N=4$, for which the curve
bends in the other direction. This is a consequence of the fact
that for $N=4$ there exists only one state (namely $\sigma_2$)
where the probability $p_B^2$ is used. This is confirmed by our
findings when we modify the definition of game $B$ such that there
is for any $N$ only one state where $p_B^2$ is used. The fact that
all curves of inversion points are symmetric upon reflection about
the plane $p_B^1=1-p_B^3$ is a consequence of the property of
Eq.~(\ref{symmetry}).

\section{Parrondo's games and the current inversion }\label{parrondo_inversion}

As stated previously, in the original Parrondo games the effect of
a current inversion when varying the mixing probability $\gamma$
is not possible. One way of understanding the reason is by means
of the relation that has been established
recently~\cite{tam03.1,tam03.2,at04.1} between the Brownian
ratchet and Parrondo's games. A fair or unfair paradoxical game
corresponds to a periodic or tilted potential respectively in the
model of a Brownian ratchet.

\begin{figure}[h]
\centerline{\includegraphics[width=13cm]{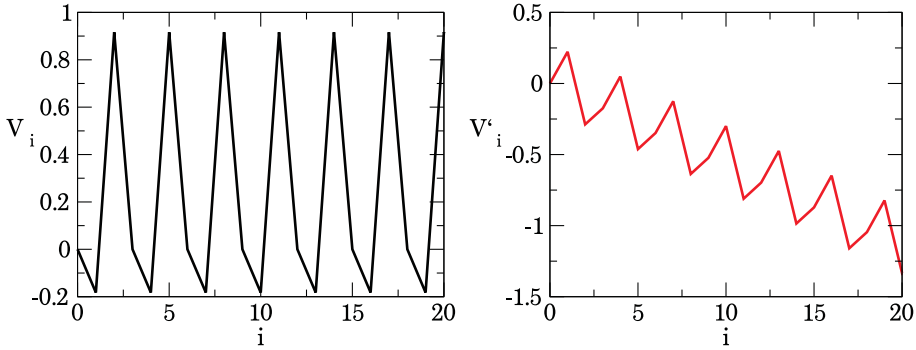}}
\caption{\label{potential_games}\textbf{\textit{a)}} Plot of the
potential related to the original Parrondo game B obtained with
the relation described in \cite{tam03.1}. \textbf{\textit{b)}}
Effective potential that we obtain when alternating between the
original Parrondo games A and B with probability
$\gamma=\frac{1}{2}$.}
\end{figure}

As an illustration, we have depicted the potential corresponding
to the original game B in Fig.~\ref{potential_games}.a. If we now
combine game B with A --which would have an associated flat
potential-- with a certain probability $\gamma$, the potential
obtained is no longer periodic, i.e., it is tilted to the right in
agreement with the direction of the flux, see
Fig.~\ref{potential_games}.b. Therefore, the question now reduces
to explain why there is no current inversion in the flashing ratchet model
when varying the rate of alternation between the potentials.

\begin{figure}[h]
\centerline{\includegraphics[width=7.5cm]{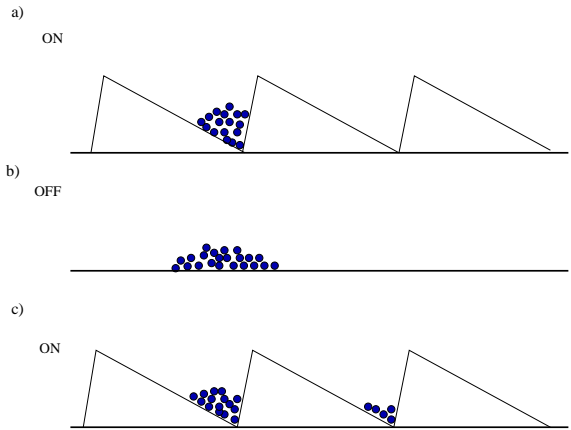}}
\caption{\label{ratchet_mechanism} Different stages of the
mechanism of rectification when switching a ratchet potential on
and off.}
\end{figure}

In the flashing ratchet model, the appearance of a flux when
alternating between a flat and an asymmetric potential is due to
a rectification process: if we consider a bunch of Brownian
particles subjected to a ratchet potential, they will tend to
remain in the potential well for a sufficiently small temperature
--see Fig.~\ref{ratchet_mechanism}.a~. When the potential is
switched off the particles start diffusing, and if we wait for
long enough, it is more likely that a small fraction of particles
will reach the vicinity of the potential well located to the
right, rather than the one on the left, due to the geometry of the
potential. Switching on and off the potential many times,
generates a net flux to the right. Whence, whatever the rate of
alternation between these two potentials, it will always be more
likely for a particle to move rightwards rather leftwards. Therefore it
will be impossible to obtain a current reversal by means of
varying only the flip rate of the potentials.


\section{Conclusions}\label{conclusions}
We have presented a new type of collective Parrondo games. These
games present, besides the Parrondo effect, a current inversion
when varying the alternation probability $\gamma$ between the two
games A and B. This phenomenon is new in the literature on
paradoxical games and the related ratchet models. Analytical
expressions for the games have been obtained for a finite number
of players using discrete--time Markov chain techniques. We have
also explained qualitatively the reason of this current inversion.

It remains as an open question the possible implications of these findings in
the field of the Brownian ratchet, as well as the possibility of finding a physical model
equivalent to this collective game.

\vspace{0.25cm}

{\bf Acknowledgments:} The authors thank C. Van den Broeck for
very helpful discussions. PA and RT acknowledge financial support
from the Ministerio de Educaci\'on y Ciencia (Spain), FEDER
projects FIS2004-5073, FIS2004-953. PA is supported from the local
governement of the Balearic Islands.

\bibliographystyle{elsart-num}
\bibliography{Bibliography}

\end{document}